\documentclass{article}
\usepackage{amsmath}
\usepackage{graphicx}
\usepackage{color}
\usepackage{verbatim}
\usepackage{footmisc}

\newcommand{\hth}{\widehat{\theta}}
\newcommand{\SE}{\mbox{SE}}

\begin{document}

\thispagestyle{empty}

\begin{center}

\textbf{\LARGE A five-decision testing procedure to infer on unidimensional parameter}

\medskip
\medskip

Aaron McDaid\footnote[1]{AMD and ZK were supported by SystemsX.ch (51RTP0 151019) and by the Swiss National Science Foundation (31003A-143914).}, Zolt\'an Kutalik$^*$, Valentin Rousson

\medskip
\medskip

Division of Biostatistics, Institute for Social and Preventive Medicine, University Hospital Lausanne, Route de la Corniche 10, 1010 Lausanne, Switzerland

\medskip
\medskip

{\tt aaron.mcdaid@gmail.com}; {\tt zoltan.kutalik@gmail.com}; {\tt valentin.rousson@chuv.ch}

\end{center}

\medskip
\medskip

\noindent \textit{Abstract}: A statistical test can be seen as a procedure to produce a decision based on observed data, where some decisions consist of rejecting a hypothesis (yielding a significant result) and some do not, and where one controls the probability to make a wrong rejection at some pre-specified significance level. Whereas traditional hypothesis testing involves only two possible decisions (to reject or not a null hypothesis), Kaiser's directional two-sided test as well as the more recently introduced Jones and Tukey's testing procedure involve three possible decisions to infer on unidimensional parameter. The latter procedure assumes that a point null hypothesis is impossible (e.g.\ that two treatments cannot have exactly the same effect), allowing a gain of statistical power. There are however situations where a point hypothesis is indeed plausible, for example when considering hypotheses derived from Einstein's theories. In this article, we introduce a five-decision rule testing procedure, which combines the advantages of the testing procedures of Kaiser (no assumption on a point hypothesis being impossible) and of Jones and Tukey (higher power), allowing for a non-negligible (typically 20\%) reduction of the sample size needed to reach a given statistical power to get a significant result, compared to the traditional approach.

\medskip


\newpage

\section{Introduction}

A statistical test can be seen as a procedure to produce a decision based on observed data (Kaiser, 1960). For example, traditional one-sided and two-sided tests to make inference on a unidimensional parameter are testing procedures with two possible decisions (to reject or not a null hypothesis). 
On the other hand, Kaiser (1960) and Jones and Tukey (2000) introduced testing procedures with three possible decisions. In this article, we propose a testing procedure with five possible decisions. 
In all these testing procedures, some decisions consist of rejecting a hypothesis, yielding a ``significant result'', and some do not. In what follows, a testing procedure is said to be valid if it allows to control the probability to make a wrong rejection, i.e.\ to reject a hypothesis which is true, the ``significance level'' $\alpha$ of a testing procedure being defined as the {\it maximal probability to make a wrong rejection}, typically set at 5\%. 

In what follows, we shall consider some unidimensional parameter $\theta$ and some 
reference value of interest $\theta_0$. In a one-sided test to the left, one attempts to reject $\theta\ge\theta_0$. In a one-sided test to the right, one attempts to reject $\theta\le\theta_0$. In a traditional two-sided test, one attempts to reject $\theta=\theta_0$, with no information whether $\theta\ge\theta_0$ or $\theta\le\theta_0$ is rejected in case of a significant result. As noted by Kaiser (1960), ``it seems difficult to imagine a problem for which this traditional test could give results of interest'' and ``to find a `significant' effect and not be able to decide in which direction this difference or effect lies, seems a sterile way to do business''. This is why he proposed instead a ``directional two-sided test'' which is equivalent to perform two one-sided tests, one to the left and one to the right, where one has the possibility to reject either $\theta\ge\theta_0$ or $\theta\le\theta_0$, depending on which one-sided test is significant. To maintain the probability of a wrong rejection at pre-specified significance level $\alpha$ in case $\theta=\theta_0$ is true, the two one-sided tests are run at the nominal significance level  $\alpha/2$. 

Most practitioners are actually applying (sometimes implicitly) a directional two-sided test when inferring on a unidimensional parameter. Some authors have objected, however, that a point hypothesis $\theta=\theta_0$ (contrary to a composite hypothesis) is almost certainly false. For example, a null hypothesis stating that the effects of two treatments A and B are equal is (in a strict sense) false since one of the two treatments A or B is inevitably superior to the other, even if not in a clinical relevant way. Jones and Tukey (2000) referred therefore to   ``the fiction of the null hypothesis'' and concluded that ``point hypotheses, while mathematically convenient, are never fulfilled in practice''. Other quotations from the literature include e.g.\ ``the null hypothesis is quasi-always false'' (Meehl, 1978), ``all we know about the world teaches us that the effects of A and B are always different - in some decimal places - for any A and B'' (Tukey, 1991), or ``in most comparative clinical trials, the point null hypothesis of no difference is not really believable'' (Freedman, 2008). Considering the point hypothesis $\theta=\theta_0$ to be impossible implies that the two one-sided tests performed in a directional two-sided test can actually be run at the nominal significance level $\alpha$ (instead of $\alpha/2$), yielding the three-decision testing procedure of Jones and Tukey (2000). 

Of course, not having to divide the nominal significance level by two when running the two one-sided tests implies a higher probability of getting a significant result, such that Jones and Tukey's testing procedure is more powerful than Kaiser's one. The price to pay for this gain of power is to assume that the point hypothesis $\theta=\theta_0$ is impossible, which might be regarded as an infinitely mild assumption. There are however situations where a point hypothesis is indeed plausible, for example in mathematics, if one considers e.g.\ the hypothesis that there are exactly 50\% of odd digits among the decimals of $\pi$,  or in particle physics, when considering the well-accepted hypothesis that an anti-electron and an electron have the same mass, or hypotheses derived from Einstein's theories, among others. In the present paper, we propose a five-decision testing procedure which combines the advantages of the procedures of Kaiser (no assumption on the point hypothesis being impossible) and of Jones and Tukey (higher power). Our five-decision rule simplifies to that of Jones and Tukey if one assumes that a point hypothesis is impossible, and even without this assumption, yields an increase of the statistical power to get a significant result compared to Kaiser's approach. Our five-decision testing procedure is described and its validity is established in Section 2. An illustration is provided in Section 3. Statistical power and sample size calculation are examined in Section 4. Section 5 contains some concluding remarks.

\section{A five-decision testing procedure}

As done in Section 1, we consider some unidimensional parameter $\theta$, e.g.\ a mean difference or a correlation, and some reference value of interest $\theta_0$ for this parameter, e.g.\ the value 0. We consider the hypotheses, $H_1:\theta\ge\theta_0$, $H_2:\theta>\theta_0$, $H_3:\theta=\theta_0$, $H_4:\theta<\theta_0$ and $H_5:\theta\le\theta_0$. While $H_1$, $H_2$, $H_4$ and $H_5$ are composite hypotheses that we shall try to reject using our testing procedure, $H_3$ is a point hypothesis which we refer to as ``the null hypothesis'', although it will be only a ``working hypothesis'' in what follows. We then consider a test statistic $T_{stat}$, a random variable with c.d.f.\ $F_{\theta}(t)=\Pr_{\theta}\{T_{stat}\le t\}$ which depends on the true value of $\theta$, and we denote by $t_{stat}$ its realization calculated from a sample of data. We make the following assumptions:

\begin{description}

\item[(A1)] the distribution of $T_{stat}$ under the null hypothesis (in what follows, the null distribution), and hence $F_{\theta_0}(t)$ is known, and let $q_{\alpha}=F^{-1}_{\theta_0}(\alpha)$ (where $0<\alpha<1$)

\item[(A2)]  $F_{\theta}(t)$ is monotone in $\theta$, such that $\theta_1<\theta_2$ implies $F_{\theta_1}(t)\ge F_{\theta_2}(t)$ (for all $t$) 

\item[(A3)] to avoid unnecessary complications in our exposition below, we consider that the null distribution is truly continuous, such that $\Pr_{\theta_0}\{T_{stat}<t\}=\Pr_{\theta_0}\{T_{stat}\le t\}$ (whatever $t$) and such that $F_{\theta_0}(q_{\alpha})=\alpha$ (for all $\alpha$).

\end{description}

\noindent Note that Assumption (A1) is needed in any statistical test involving a point null hypothesis, assumption (A2) is classical (ensuring e.g.\ unbiased tests and monotonicity of statistical power), whereas assumption (A3) could be relaxed (although this would require more complicated notations). Given a pre-specified significance level $\alpha$, our five-decision testing procedure is then defined as follows: 

\medskip

\begin{center}
\begin{tabular}{ccc}
\hline
decision & event & hypothesis rejected   \\
\hline
1 & $t_{stat}<q_{\alpha/2}$ & $H_1:\theta\ge\theta_0$  \\
2 & $q_{\alpha/2}\le t_{stat}<q_{\alpha}$ & $H_2:\theta>\theta_0$ \\
3 & $q_{\alpha}\le t_{stat}\le q_{1-\alpha}$ & none  \\
4 & $q_{1-\alpha}<t_{stat}\le q_{1-\alpha/2}$ & $H_4:\theta<\theta_0$ \\
5 & $q_{1-\alpha/2}<t_{stat}$ & $H_5:\theta\le\theta_0$ \\
\hline
\end{tabular}
\end{center}

\medskip 

\noindent To ensure mutually exclusive decisions, we consider $0<\alpha \le 0.5$. Note that the first, second, fourth and fifth decisions result in the rejection of a hypothesis whereas the third decision does not. Note also that some rejections are stronger than other, rejection of $H_1$ implying rejection of $H_2$ and rejection of $H_5$ implying rejection of $H_4$.
As in Kaiser (1960) and Jones and Tukey (2000), when a hypothesis is rejected we consider the complementary hypothesis to be implicitly accepted (which is particularly simple to define since the rejected hypothesis involves only a unidimensional parameter). Thus, rejection of $H_1$, $H_2$, $H_4$ or $H_5$ implicitly implies acceptance of $H_4$, $H_5$, $H_1$ and $H_2$ respectively.

Of note, our five-decision testing procedure could be formulated as a combination of three traditional tests: two one sided-tests, one to the left (OSL) where one tries to reject $H_1:\theta\ge\theta_0$, one to the right (OSR) where one tries to reject $H_5:\theta\le\theta_0$, and one traditional two-sided test (TS), where one tries to reject the null hypothesis $H_3:\theta=\theta_0$, each of them run at the same significance level $0<\alpha\le0.5$, with five possibilities corresponding to our five decisions as follows:

\bigskip

\begin{center}
\begin{tabular}{ccc}
\hline
decision & outcome of traditional tests & hypothesis rejected  in the \\
&& five-decision testing procedure \\
\hline
1 & reject both $H_1$ (with OSL) and $H_3$ (with TS) & $H_1:\theta\ge\theta_0$  \\
2 & reject $H_1$ (with OSL), not $H_3$ (with TS) & $H_2:\theta>\theta_0$ \\
3 & reject neither $H_1$, $H_5$ nor $H_3$ (with OSL, OSR and TS)  & none  \\
4 & reject $H_5$ (with OSR), not $H_3$ (with TS) & $H_4:\theta<\theta_0$ \\
5 & reject both $H_5$ (with OSR) and $H_3$ (with TS) & $H_5:\theta\le\theta_0$ \\
\hline
\end{tabular}
\end{center}

\bigskip 

Although in general a testing procedure obtained as a combination of tests that control the type I error at some level $\alpha$ is not guaranteed to control type I error at $\alpha$, we demonstrate below that our five-decision testing procedure is valid. Recall that a testing procedure is valid if the probability to make a wrong rejection cannot exceed $\alpha$, whatever the true value of $\theta$. Note first that if $\theta=\theta_0$, one gets a wrong rejection when the first or fifth decision occurs, i.e.\ either when $t_{stat}<q_{\alpha/2}$ or when $t_{stat}>q_{1-\alpha/2}$. This kind of wrong rejection is known as a type I error. The probability that this happens is given by:
\[
\mbox{$\Pr_{\theta_0}$}\{t_{stat}< q_{\alpha/2}\}+\mbox{$\Pr_{\theta_0}$}\{t_{stat}> q_{1-\alpha/2}\}=
F_{\theta_0}(q_{\alpha/2})+1-F_{\theta_0}(q_{1-\alpha/2})=\alpha/2+1-(1-\alpha/2)=\alpha.
\]
If $\theta<\theta_0$, one gets a wrong rejection when the fourth or fifth decision occurs, i.e.\ when $t_{stat}>q_{1-\alpha}$.  In case of decision 5, such a wrong rejection is sometimes referred to as a type III error (Kimball, 1957; Leventhal and Huynh, 1996; Shaffer, 2002) or ``operational type I error'' (Senn, 2007, p.\ 188). The probability that this happens is given by:
\[
\mbox{$\Pr_{\theta}$}\{t_{stat}> q_{1-\alpha}\}=
1-F_{\theta}(q_{1-\alpha})\le 1-F_{\theta_0}(q_{1-\alpha})=1-(1-\alpha)=\alpha. 
\]
Thus, this probability is unknown but not larger than $\alpha$. If $\theta>\theta_0$, one gets a wrong rejection when the first or second decision occurs, i.e.\ when $t_{stat}<q_{\alpha}$.  In case of decision 1, this is another example of type III error. The probability that this happens is given by:
\[
\mbox{$\Pr_{\theta}$}\{t_{stat}<q_{\alpha}\}=
F_{\theta}(q_{\alpha})\le F_{\theta_0}(q_{\alpha})=\alpha. 
\]
Thus, whatever the true value of $\theta$, the probability to make a wrong rejection is unknown, but bounded by $\alpha$, ensuring the validity of the testing procedure, which is either correctly sized if $\theta=\theta_0$, or conservative if $\theta<\theta_0$ or $\theta>\theta_0$. If the null distribution is known only approximately, the testing procedure is still approximately valid. 

A typical example where Assumptions (A1)--(A3) hold and where the null distribution is known exactly (under some conditions such as normality and homoscedasticity) are $t$-tests. Another example where the null distribution is known approximately are Wald tests. In that case, the test statistic is defined as $T_{stat}=(\hth-\theta_0)/\SE(\hth)$, where $\hth$ is a consistent and (asymptotically) normally distributed estimate of $\theta$ and $\SE(\hth)$ is (a consistent estimate of) the standard error of $\hth$ (which is ideally calculated under the null hypothesis), such that the null distribution is approximately standard normal. One has in that case $q_{\alpha}\approx z_{\alpha}$, where $z_{\alpha}=\Phi^{-1}(\alpha)$ and 
$\Phi(t)$ refers to the c.d.f. of a standard normal distribution. With $\alpha=5\%$, decisions 1--5 are taken when respectively $t_{stat}<-1.96$, $-1.96\le t_{stat}<-1.645$, $-1.645\le t_{stat}\le 1.645$, $1.645<t_{stat}\le 1.96$ and $1.96<t_{stat}$.

\section{Illustration}

To illustrate our testing procedure, we consider the ``ChickWeight'' data set available in the R base package. In that data set, 50 chicks have been followed up during the first three weeks of life. The chicks received different experimental protein diets (20 received diet 1, 10 diet 2, 10 diet 3 and 10 diet 4) and have been weighed every two days. Figure 1 shows boxplots of the weight measured after 20 days for the 10 chicks which received diet 2, and for the 10 chicks which received diet 3. One can see that the sample mean was higher with diet 3 than with diet 2 (258.9 vs 205.6 grams), whereas sample standard deviations were similar (65.2 vs 70.3 grams), a pooled standard deviation being obtained as $((65.2^2+70.3^2)/2)^{1/2}=67.8$. Let $\mu_2$ and $\mu_3$ denote the true means in these two groups and let $\theta=\mu_3-\mu_2$. The test statistic of a two-sample $t$-test to try to reject the equality of the two means is given by $t_{stat}=(10/2)^{1/2}\cdot(258.9-205.6)/67.8=1.76$, yielding a (two-sided) $p$-value of $p=0.096$. Recall that the null distribution is here a $t$-distribution with 18 degrees of freedom, for which the quantile 97.5\% is given by $q_{0.975}=2.10$ and the quantile 95\% by $q_{0.95}=1.73$. In what follows, we discuss the results of the different testing procedures presented above, run each at the $\alpha=5\%$ significance level. 

\begin{figure}[t]
\begin{center}
\hspace*{-1cm}\includegraphics[scale=1]{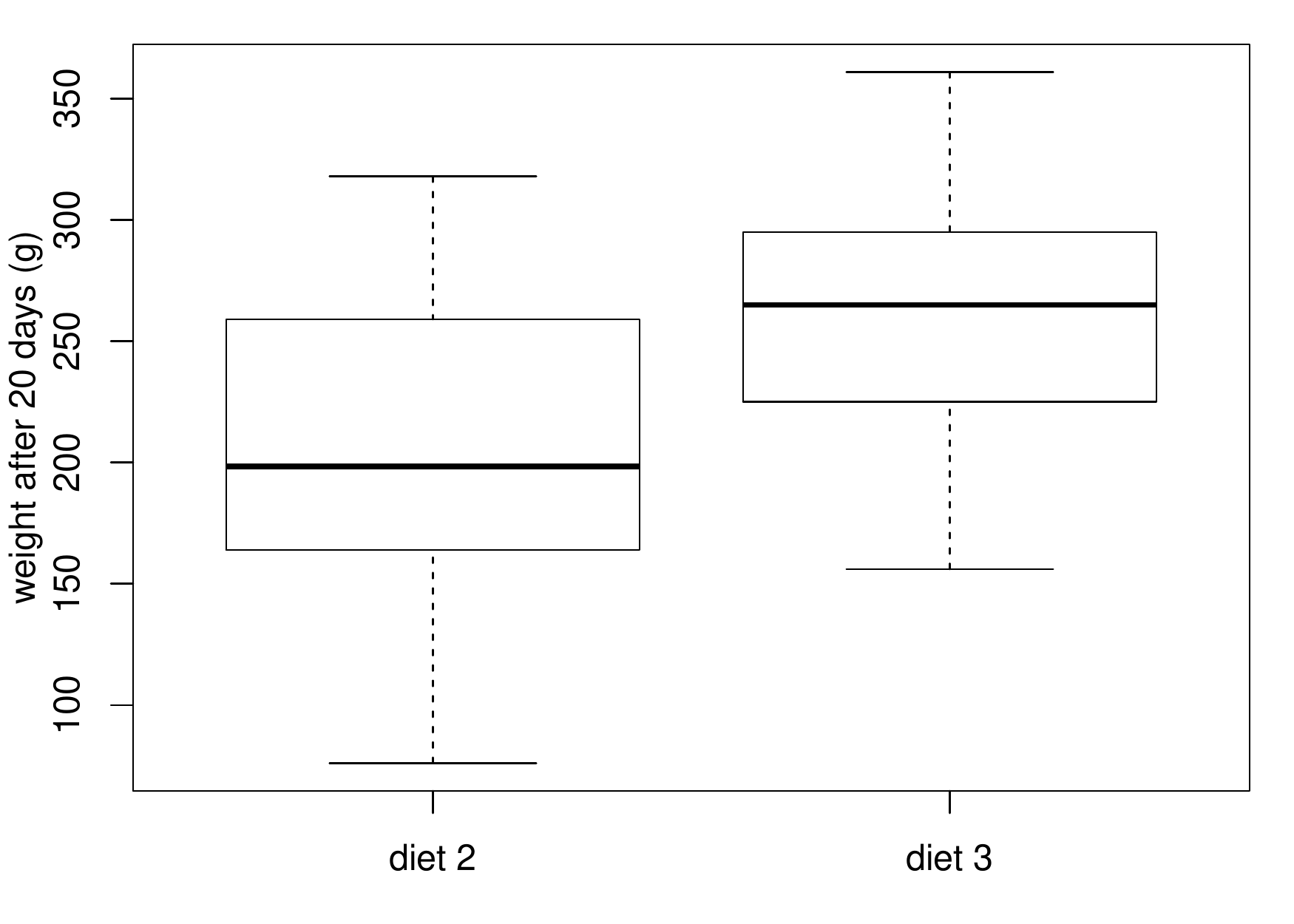}
\end{center}
\caption{\textit{Boxplots of the weight after 20 days for chicks nourished with diet 2 and diet 3. }}
\end{figure}

\begin{itemize}

\item {\bf Kaiser's testing procedure} (directional two sided-test): Since $|t_{stat}|\le q_{0.975}$, we are not able to reject the null hypothesis $H_3:\theta=0$ of no mean difference between the two groups (as already indicated by the fact that $0.05\le p$). 

\item {\bf Jones and Tukey's testing procedure} (two one-sided tests): Since $q_{0.95}<t_{stat}$, and assuming that the null hypothesis $H_3:\theta=0$ is impossible, we are able to reject the hypothesis $H_5:\theta\le0$ (and thus to conclude that the true mean with diet 3 is strictly higher than the true mean with diet 2). One could have reached the same conclusion without looking at the critical value $q_{0.95}$ by noting that $p<0.10$ and that the sample mean is higher with diet 3 than with diet 2.

\item {\bf Five-decision testing procedure}: Since $q_{0.95}<t_{stat}\le q_{0.975}$, and without assuming that the null hypothesis is impossible, we are able to reject the hypothesis $H_4:\theta<0$ (and thus to conclude that the true mean with diet 3 is at least as high than the true mean with diet 2). One could have reached the same conclusion without looking at the critical values $q_{0.95}$ and $q_{0.975}$ by noting that $0.05\le p<0.10$ and that the sample mean is higher with diet 3 than with diet 2. On the other hand, if we assume that the null hypothesis is impossible, rejecting $H_4$ is equivalent to rejecting $H_5$ such that one gets the same conclusions as with Jones and Tukey's testing procedure. 

\end{itemize}

Figure 2 illustrates what happens when performing the five-decision testing procedure at different significance levels, showing the decision achieved in function of the value of $t_{stat}$ in the context of our example (i.e.\ a two-sample $t$-test where the null distribution is a $t$-distribution with 18 degrees of freedom). With $t_{stat}=1.76$, while we just saw that one rejects $H_4:\theta<0$ at the 5\% significance level, one can see on that figure that one rejects $H_5:\theta\le 0$ at the 10\% significance level (which is a stronger rejection than to reject $H_4$), whereas no hypothesis can be rejected at the 1\% significance level.   

\begin{figure}[tbph]
\begin{center}
\hspace*{-1cm}\includegraphics[scale=1]{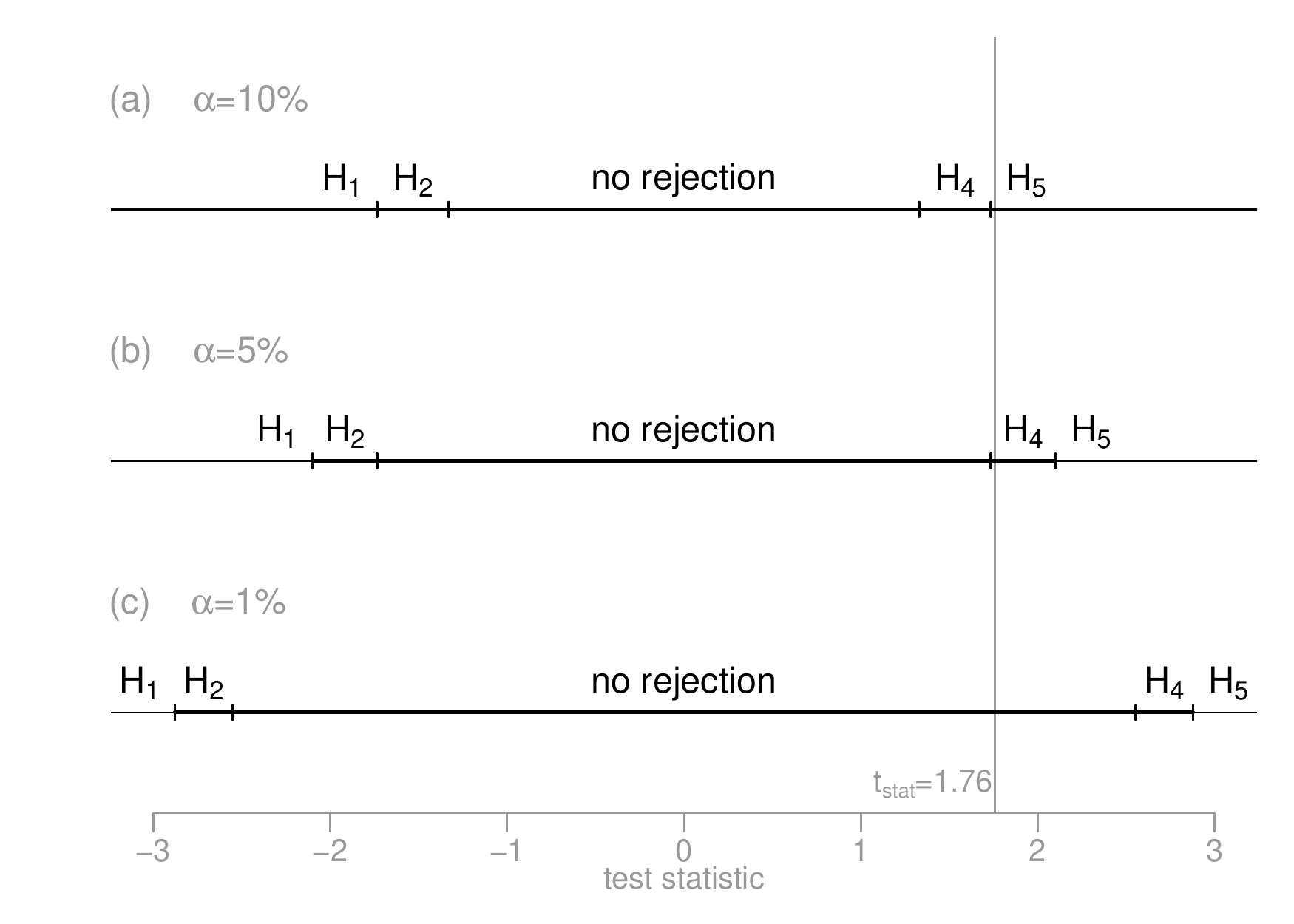}
\end{center}
\caption{\textit{
Decision (rejection) achieved when using the five-decision testing procedure at significance level (a) $\alpha=10\%$, (b) $\alpha=5\%$, and (c) $\alpha=1\%$, depending on the the value of the test statistic, in the context of our example (i.e.\ a two-sample $t$-test where the null distribution is a $t$-distribution with 18 degrees of freedom). With a calculated value of $t_{stat}=1.76$ (vertical line), one rejects $H_5:\theta\le0$ at the 10\% level, one  rejects $H_4:\theta<0$ at the 5\% level, while no hypothesis can be rejected at the 1\% level (where $\theta$ is the mean difference of weights between chicks nourished with diet 3 and chicks nourished with diet 2).}}
\end{figure}

It is well known that one rejects the null hypothesis $\theta=0$ in a traditional two-sided two-sample $t$-test at the 5\% significance level if and only if the value $0$ is lying outside a 95\% confidence interval for $\theta$. In Kaiser's testing procedure, one then rejects $\theta\ge0$ if $0$ is lying on the right of the confidence interval, and one rejects $\theta\le0$ if $0$ is lying on the left of the confidence interval. To get informed about the outcome of the five-decision testing procedure at the 5\% significance level, one needs to calculate a 90\%, in addition to a 95\% confidence interval for $\theta$, and to proceed as follows:

\begin{description}

\item[(decision 1)] reject $H_1: \theta\ge0$ if the value 0 is found on the right of the 95\% confidence interval,

\item[(decision 2)] reject $H_2: \theta>0$ if the value 0 is found on the right of the 90\%, but within the 95\% confidence interval,

\item[(decision 3)] no rejection if the value is 0 is found within the 90\% confidence interval,

\item[(decision 4)] reject $H_4: \theta<0$ if the value 0 is found on the left of the 90\%, but within the 95\% confidence interval,

\item[(decision 5)] reject $H_5: \theta\le0$ if the value 0 is found on the left of the 95\% confidence interval.

\end{description}

\noindent  In our example, a 95\% confidence interval for $\theta$ is given by $[-10.4;117.0]$ while a 90\% confidence interval for $\theta$ is given by $[0.7;105.9]$. 
Since the value $0$ belongs to the 95\% confidence interval while being on the left of the 90\% confidence interval, one (again) rejects $H_4:\theta<0$ at the 5\% significance level.

\section{Statistical power and sample size calculation}

Whereas the significance level $\alpha$ of a testing procedure is the (maximal) probability to (wrongly) reject a hypothesis which is true, the statistical power $\psi$ can be defined as the probability to (correctly) reject a hypothesis which is false. We provide below formulas for the statistical power achieved with our five-decision testing procedure in the case of a Wald test. Note that formulae provided in this section are asymptotical and may require large sample sizes to become accurate.

\begin{itemize}

\item in the case $\theta<\theta_0$, the probability to (correctly) reject $H_1:\theta\ge\theta_0$ (decision 1) is given by:
\[
\psi_1=\mbox{$\Pr_{\theta}$}\{t_{stat}<z_{\alpha/2}\}=\Phi(z_{\alpha/2}+(\theta_0-\theta)/\SE(\hth))
\]

\item in the case $\theta\le\theta_0$, the probability to (correctly) reject $H_2:\theta>\theta_0$ (decision 1 or 2) is given by:
\[
\psi_2=\mbox{$\Pr_{\theta}$}\{t_{stat}<z_{\alpha}\}=\Phi(z_{\alpha}+(\theta_0-\theta)/\SE(\hth))
\]

\item in the case $\theta\ge\theta_0$, the probability to (correctly) reject $H_4: \theta<\theta_0$ (decision 4 or 5) is given by:
\[
\psi_4=\mbox{$\Pr_{\theta}$}\{t_{stat}>z_{\alpha}\}=\Phi(z_{\alpha}+(\theta-\theta_0)/\SE(\hth))
\]

\item in the case $\theta>\theta_0$, the probability to (correctly) reject $H_5: \theta\le\theta_0$ (decision 5) is given by:
\[
\psi_5=\mbox{$\Pr_{\theta}$}\{t_{stat}>z_{\alpha/2}\}=\Phi(z_{\alpha/2}+(\theta-\theta_0)/\SE(\hth)).
\]

\end{itemize}

\noindent As an example, let consider $\alpha=5\%$ and a case with $\theta>\theta_0$, where the difference between $\theta$ and $\theta_0$ expressed in ``standard error units'' is given by $(\theta-\theta_0)/\SE(\hth)=2.5$. The probability to (correctly) reject $H_5:\theta\le\theta_0$ is given by $\psi_5=\Phi(-1.96+2.5)=70.5\%$, whereas the probability to (correctly) reject $H_4:\theta<\theta_0$ is given by $\psi_4=\Phi(-1.645+2.5)=80.4\%$. The power to reject a strict inequality ($H_4$) is logically higher than the power to reject a non-strict inequality ($H_5$). Note also that if one considers the null hypothesis to be impossible, rejecting $H_4$ will be equivalent to rejecting $H_5$, enabling an increase of statistical power from 70.5\% to 80.4\%.

An increase of statistical power allows in turn a reduction of the sample size $n$ needed to reach a given statistical power when designing a study. In the case of a Wald test where the standard error of the estimate $\hth$ is given by $\SE(\hth)=\tau/\sqrt{n}$ (not depending on the true value of $\theta$), the sample size needed to reject a non-strict inequality ($H_1:\theta\ge\theta_0$ or $H_5:\theta\le\theta_0$) with given probability (power) $\psi$ is given by:
\begin{equation}
n=\frac{(z_{1-\alpha/2}+z_{\psi})^2\tau^2}{(\theta-\theta_0)^2}
\label{nlarge}
\end{equation}
whereas  the sample size needed to reject a strict inequality ($H_2:\theta>\theta_0$ or $H_4:\theta<\theta_0$) with given probability (power) $\psi$ is given by: 
\begin{equation}
n=\frac{(z_{1-\alpha}+z_{\psi})^2\tau^2}{(\theta-\theta_0)^2}.
\label{nstrict}
\end{equation}
Compared to a traditional sample size calculation \eqref{nlarge} for rejecting only a non-strict inequality, the sample size calculated via \eqref{nstrict} for rejecting a strict inequality, as enabled in our testing procedure, yields a relative reduction of sample size of:
\begin{equation}
\frac{(z_{1-\alpha/2}+z_{\psi})^2-(z_{1-\alpha}+z_{\psi})^2}{(z_{1-\alpha/2}+z_{\psi})^2}.
\label{ngain}
\end{equation}
The following table provides such examples of sample size reduction \eqref{ngain} in function of $\psi$ and $\alpha$: 

\medskip

\begin{center}
\begin{tabular}{c|ccccc}
\hline
& $\psi=50\%$ & $\psi=80\%$ & $\psi=90\%$ & $\psi=95\%$ & $\psi=99\%$ \\
\hline
$\alpha=5\%$ & 30\% & 21\% & 18\% & 17\% & 14\% \\
$\alpha=1\%$ & 18\% & 14\% & 13\% & 11\% & 10\% \\
$\alpha=0.5\%$ & 16\% & 12\% & 11\% & 10\% & 9\% \\
$\alpha=0.1\%$ & 12\% & 9\% & 9\% & 8\% & 7\% \\
\hline
\end{tabular}
\end{center}

\medskip

\noindent  
Thus, having settled for rejecting a strict rather than a non-strict inequality, which will be of particular interest for those assuming the null hypothesis to be impossible (since both rejections are then equivalent), enables e.g.\ a reduction of sample size of 21\% in a study targeting a statistical power of $\psi=80\%$ with $\alpha=5\%$. 

To further illustrate such sample size reduction and that formulas \eqref{nlarge} and \eqref{nstrict} can also be useful when using an exact test, we consider an example where one would attempt to show that a treatment A is superior to a treatment B via a two-sample $t$-test. The parameter of interest is here $\theta=\mu_A-\mu_B$, where $\mu_A$ and $\mu_B$ represent the true means of some continuous health outcome, characterizing the effects of treatments A and B, the reference value being $\theta_0=0$ and the null hypothesis $H_3:\mu_A=\mu_B$. The test statistic is given by $t_{stat}=\sqrt{n/2}(\bar{x}_A-\bar{x}_B)/s$ with $s^2=(s_A^2+s_B^2)/2$, where $\bar{x}_A$ and $\bar{x}_B$ denote the sample means and $s_A^2$ and $s_B^2$ the sample variances of the health outcome calculated from two samples of size $n$. Assuming a normal distribution and a same variance $\sigma^2$ for both treatments, the null distribution is a $t$-distribution with $2n-2$ degrees of freedom, which can be approximated by a standard normal distribution for a large $n$, as in a Wald test. If the goal is to reject $H_5:\mu_A\le \mu_B$, expecting a ``medium'' treatment effect expressed as $(\mu_A-\mu_B)/\sigma=0.5$ (Cohen, 1988), using a significant level $\alpha=5\%$, and targeting a statistical power of $\psi=80\%$, one may calculate a sample size via \eqref{nlarge} of (noting that $SE(\hth)=\sqrt{2\sigma^2/n}$, such that $\tau^2=2\sigma^2$):
\[
n=\frac{(z_{1-\alpha/2}+z_{\psi})^2\cdot 2\sigma^2}{(\mu_A-\mu_B)^2}
=\frac{2(1.96+0.84)^2}{0.5^2}=63.
\]
Now, if one considers that the null hypothesis $H_3:\mu_A=\mu_B$ is impossible (i.e.\ that the two treatments cannot have exactly the same effect), one will content to reject the hypothesis $H_4:\mu_A<\mu_B$ and one will calculate a sample size via \eqref{nstrict} of: 
\[
n=\frac{(z_{1-\alpha}+z_{\psi})^2\cdot2\sigma^2}{(\mu_A-\mu_B)^2}
=\frac{2(1.645+0.84)^2}{0.5^2}=50
\]
achieving a $(63-50)/63=21\%$ reduction of sample size. Under such assumptions, one can check via simulation that the probability to reject $H_5:\mu_A\le\mu_B$ (decision 5 from our testing procedure) via a two-sample $t$-test with two groups of size $n=63$, i.e.\ to get $t_{stat}>1.979$ (the quantile 97.5\% of a $t$-distribution with 124 degrees of freedom), is about 79.3\%, whereas the probability to reject $H_4:\mu_A<\mu_B$ (decision 4 or decision 5 from our testing procedure) with two groups of size $n=50$, i.e.\ to get $t_{stat}>1.661$ (the quantile 95\% of a $t$-distribution with 98 degrees of freedom) is about 79.7\% (estimated from 100'000 simulations), both pretty close to the targeted power of 80\%.

\section{Conclusions}

The formulation of hypothesis testing proposed by Kaiser (1960) and Jones and Tukey (2000) allows one to perform on the same data two traditional one-sided tests, trying to reject two different composite hypotheses. Since the corresponding alternative hypotheses are also one-sided, interpretation of a significant result is straightforward. Not having to divide by two the nominal significance level when performing these two tests, as advocated by Jones and Tukey (2000), also implies  ``the abolition, once and for all, of the controversy over whether a one-sided or two-sided test is appropriate'' (Freedman, 2008, 2009). Note that a similar procedure is used in bioequivalence studies, where two one-sided tests, run to reject two disjoint composite hypotheses, are also typically conducted at 5\%, and where the calculation of a 90\% (instead of a 95\%) confidence interval has been advocated (Schuirmann, 1987; Westlake, 1981).

In this paper, we have introduced a five-decision testing procedure which can be seen as an extension of both approaches. On the one hand, this is an extension of Kaiser's testing procedure from a three- to a five-decision testing procedure (our decisions 2 and 4 being absent from Kaiser's procedure, in fact merged with our decision 3). On the other hand, our five-decision testing procedure reduced to Jones and Tukey's three-decision testing procedure if the null hypothesis is considered to be impossible (decision 1 being then equivalent to decision 2, and decision 5 being equivalent to decision 4). Importantly, the five-decision testing procedure can be used both by those who believe in the plausibility of the null hypothesis and those who do not. For the former, our approach is  still more powerful than Kaiser's one, allowing to distinguish between the rejection of a strict and of a non-strict inequality. For the latter, it is as powerful as Jones and Tukey's approach, allowing a non negligible reduction of the sample size needed to reach a given statistical power compared to a traditional sample size calculation, e.g.\ of 21\%, as illustrated in our example of Section 4, although this calculation was based on asymptotic formulae and may need large sample sizes to become effective. 

Although our approach is clearly frequentist, it is interesting to see how allowing the option of believing or not in the plausibility of a null hypothesis reflects Bayesian thinking. While controlling the probability to make a wrong rejection is not a Bayesian concept, one considers in a Bayesian context a prior and a posterior distribution for a parameter $\theta$, such that it is also possible to assign prior and calculate posterior probabilities associated to the hypotheses $H_1$, $H_2$, $H_3$, $H_4$ and $H_5$ considered above. In that context, those believing in the plausibility of the null hypothesis should assign a non zero prior probability to $H_3$, translating to a non-continuous prior distribution for $\theta$ (with a point mass at the reference value $\theta_0$) in Bayesian inference. As a consequence, the posterior probabilities for (and hence Bayesian inference about) $H_1$ and $H_5$ will be different than for respectively $H_2$ and $H_4$. In the other case ($H_3$ being impossible), the posterior probabilities for (and hence Bayesian inference about) $H_1$ and $H_5$ are the same than for respectively $H_2$ and $H_4$. Therefore, a common point in Bayesian inference and our five-decision testing procedure is that believing or not in the plausibility of the null hypothesis $H_3$ affects inference on $H_1$, $H_2$, $H_4$ and $H_5$. 

\medskip


\end{document}